\documentclass[leqno,12pt]{amsart}
\usepackage{amssymb} 
\theoremstyle{plain} 
\newtheorem{theorem}{\indent\sc Theorem}[section] 

\newtheorem{proposition}[theorem]{\indent\sc Proposition}

\newtheorem{question}[theorem]{\indent\sc Question}
\theoremstyle{definition} 
\newtheorem{definition}[theorem]{\indent\sc Definition}

\begin{document}

\title[Fourth-order differential equation]{Fourth-order ordinary differential equation obtained by similarity reduction of the\\
 modifed Sawada-Kotera equation \\}
\author{Yusuke Sasano }

\renewcommand{\thefootnote}{\fnsymbol{footnote}}
\footnote[0]{2000\textit{ Mathematics Subjet Classification}.
34M55; 34M45; 58F05; 32S65.}

\maketitle

\begin{abstract}
We study a one-parameter family of the fourth-order ordinary differential equations obtained by similarity reduction of the modifed Sawada-Kotera equation. We show that the birational transformations take this equation to the polynomial Hamiltonian system in dimension four. We make this polynomial Hamiltonian from the viewpoint of accessible singularity and local index. We also give its symmetry and holomorphy conditions. These properties are new. Moreover, we introduce a symmetric form in dimension five for this Hamiltonian system by taking the two invariant divisors as the dependent variables. Thanks to the symmetric form, we show that this system admits the affine Weyl group symmetry of type $A_2^{(2)}$ as the group of its B{\"a}cklund transformations.
\end{abstract}

\section{Introduction}

In this paper, we study a 1-parameter family of the fourth-order ordinary differential equations in the polynomial class (see \cite{Cosgrove}) explicitly given by
\begin{equation}\label{eq:1}
\frac{d^4u}{dt^4}=-5\frac{du}{dt} \frac{d^2u}{dt^2}+5u^2 \frac{d^2u}{dt^2}+5u\left(\frac{du}{dt} \right)^2-u^5+tu+\alpha \quad (\alpha \in {\Bbb C}).
\end{equation}
This equation is well-known to be a group invariant reduction of the Modified Sawada-Kotera equation (see \cite{Cosgrove}):
\begin{equation}\label{SKeq}
g_{wwwww}=-5(g_w-g^2)g_{www}-5(g_{ww})^2+20g g_{w} g_{ww}+5(g_{w})^3-5g^4 g_{w}+g_{s}
\end{equation}
found by Fordy and Gibbons \cite{Fordy}, who observed that this equation is also the modified version of the Kaup-Kupershmidt equation. The hierarchy of Painlev\'e equations built on this equation has been studied recently by Joshi and Pickering.

We note that this equation appears as the equation F-XVIII in Cosgrove's classification of the fourth-order ordinary differential equations in the polynomial class having the Painlev\'e property (see \cite{Cosgrove}).

At first, we show that the birational transformations (see Section 2) take the equation \eqref{eq:1} to the polynomial Hamiltonian system in dimension four. We make this polynomial Hamiltonian from the viewpoint of accessible singularity and local index (see Section 5).

This Hamiltonian system can be considered as a 1-parameter family of coupled Painlev\'e systems in dimension four.

We also study its symmetry and holomorphy conditions. These B{\"a}cklund transformations satisfy
\begin{equation}
s_i(g)=g+\frac{\alpha_i}{f_i}\{f_i,g\}+\frac{1}{2!} \left(\frac{\alpha_i}{f_i} \right)^2 \{f_i,\{f_i,g\} \}+\cdots \quad (g \in {\Bbb C}(t)[q_1,p_1,q_2,p_2]),
\end{equation}
where poisson bracket $\{,\}$ satisfies the relations:
$$
\{p_1,q_1\}=\{p_2,q_2\}=1, \quad the \ others \ are \ 0.
$$
Since these B{\"a}cklund transformations have Lie theoretic origin, similarity reduction of a Drinfeld-Sokolov hierarchy admits such a B{\"a}cklund symmetry.

These properties of its symmetry and holomorphy conditions are new.

Moreover, we introduce a symmetric form in dimension five for this Hamiltonian system by taking the two invariant divisors as the dependent variables (see Section 7). Thanks to the symmetric form, we show that this system admits the affine Weyl group symmetry of type $A_2^{(2)}$ as the group of its B{\"a}cklund transformations.

\section{Polynomial Hamiltonian system}
\begin{theorem}\label{th1.1}
The birational transformations
\begin{equation}\label{eq:2}
  \left\{
  \begin{aligned}
   q_1 =&-u,\\
   p_1 =&-\left(\frac{d^3u}{dt^3}-\frac{1}{2}u^4-u^2 \frac{du}{dt}+\frac{3}{2}\left(\frac{du}{dt} \right)^2+2u \frac{d^2u}{dt^2}+\frac{t}{2} \right),\\
   q_2 =&-\left(\frac{d^2u}{dt^2}+2u\frac{du}{dt} \right),\\
   p_2 =&-\left(\frac{du}{dt}+u^2 \right)
   \end{aligned}
  \right. 
\end{equation}
take the equation \eqref{eq:1} to the Hamiltonian system
\begin{equation}\label{eq:3}
  \left\{
  \begin{aligned}
   \frac{dq_1}{dt} =&\frac{\partial H}{\partial p_1}=q_1^2+p_2,\\
   \frac{dp_1}{dt} =&-\frac{\partial H}{\partial q_1}=-2q_1p_1-\alpha-\frac{1}{2},\\
   \frac{dq_2}{dt} =&\frac{\partial H}{\partial p_2}=-\frac{1}{2}p_2^2+p_1+\frac{t}{2},\\
   \frac{dp_2}{dt} =&-\frac{\partial H}{\partial q_2}=q_2
   \end{aligned}
  \right. 
\end{equation}
with the polynomial Hamiltonian:
\begin{align}\label{eq:4}
\begin{split}
H=&K(q_1,p_1;\alpha)+H_I(q_2,p_2,t)+p_1p_2\\
=&q_1^2p_1+\left(\frac{1}{2}+\alpha \right)q_1-\frac{1}{6}p_2^3+\frac{t}{2}p_2-\frac{q_2^2}{2}+p_1p_2.
\end{split}
\end{align}
\end{theorem}
The symbols $K(x,y;\alpha)$ and $H_I(z,w,t)$ denote
\begin{align}\label{eq:5}
\begin{split}
K(x,y;\alpha)=&x^2y+\left(\frac{1}{2}+\alpha \right)x\\
H_I(z,w,t)=&-\frac{1}{6}w^3+\frac{t}{2}w-\frac{z^2}{2}.
\end{split}
\end{align}
The system with the Hamiltonian $K(x,y;\alpha)$ has itself as its first integral, and  $H_I(z,w,t)$ denotes the Painlev\'e I Hamiltonian.

This Hamiltonian system can be considered as a 1-parameter family of coupled Painlev\'e systems in dimension four.

Before we will prove Theorem \ref{th1.1}, we review the notion of accessible singularity and local index.

\section{Accessible singularity and local index}
Let us review the notion of {\it accessible singularity}. Let $B$ be a connected open domain in $\Bbb C$ and $\pi : {\mathcal W} \longrightarrow B$ a smooth proper holomorphic map. We assume that ${\mathcal H} \subset {\mathcal W}$ is a normal crossing divisor which is flat over $B$. Let us consider a rational vector field $\tilde v$ on $\mathcal W$ satisfying the condition
\begin{equation*}
\tilde v \in H^0({\mathcal W},\Theta_{\mathcal W}(-\log{\mathcal H})({\mathcal H})).
\end{equation*}
Fixing $t_0 \in B$ and $P \in {\mathcal W}_{t_0}$, we can take a local coordinate system $(x_1,\ldots ,x_n)$ of ${\mathcal W}_{t_0}$ centered at $P$ such that ${\mathcal H}_{\rm smooth \rm}$ can be defined by the local equation $x_1=0$.
Since $\tilde v \in H^0({\mathcal W},\Theta_{\mathcal W}(-\log{\mathcal H})({\mathcal H}))$, we can write down the vector field $\tilde v$ near $P=(0,\ldots ,0,t_0)$ as follows:
\begin{equation*}
\tilde v= \frac{\partial}{\partial t}+g_1 
\frac{\partial}{\partial x_1}+\frac{g_2}{x_1} 
\frac{\partial}{\partial x_2}+\cdots+\frac{g_n}{x_1} 
\frac{\partial}{\partial x_n}.
\end{equation*}
This vector field defines the following system of differential equations
\begin{equation}\label{39}
\frac{dx_1}{dt}=g_1(x_1,\ldots,x_n,t),\ \frac{dx_2}{dt}=\frac{g_2(x_1,\ldots,x_n,t)}{x_1},\cdots, \frac{dx_n}{dt}=\frac{g_n(x_1,\ldots,x_n,t)}{x_1}.
\end{equation}
Here $g_i(x_1,\ldots,x_n,t), \ i=1,2,\dots ,n,$ are holomorphic functions defined near $P=(0,\dots ,0,t_0).$

\begin{definition}\label{Def1}
With the above notation, assume that the rational vector field $\tilde v$ on $\mathcal W$ satisfies the condition
$$
(A) \quad \tilde v \in H^0({\mathcal W},\Theta_{\mathcal W}(-\log{\mathcal H})({\mathcal H})).
$$
We say that $\tilde v$ has an {\it accessible singularity} at $P=(0,\dots ,0,t_0)$ if
$$
x_1=0 \ {\rm and \rm} \ g_i(0,\ldots,0,t_0)=0 \ {\rm for \rm} \ {\rm every \rm} \ i, \ 2 \leq i \leq n.
$$
\end{definition}

If $P \in {\mathcal H}_{{\rm smooth \rm}}$ is not an accessible singularity, all solutions of the ordinary differential equation passing through $P$ are vertical solutions, that is, the solutions are contained in the fiber ${\mathcal W}_{t_0}$ over $t=t_0$. If $P \in {\mathcal H}_{\rm smooth \rm}$ is an accessible singularity, there may be a solution of \eqref{39} which passes through $P$ and goes into the interior ${\mathcal W}-{\mathcal H}$ of ${\mathcal W}$.

Here we review the notion of {\it local index}. Let $v$ be an algebraic vector field with an accessible singular point $\overrightarrow{p}=(0,\ldots,0)$ and $(x_1,\ldots,x_n)$ be a coordinate system in a neighborhood centered at $\overrightarrow{p}$. Assume that the system associated with $v$ near $\overrightarrow{p}$ can be written as
\begin{align}\label{b}
\begin{split}
\frac{d}{dt}\begin{pmatrix}
             x_1 \\
             x_2 \\
             \vdots\\
             x_{n-1} \\
             x_n
             \end{pmatrix}=\frac{1}{x_1}\left\{\begin{bmatrix}
             a_{11} & 0 & 0 & \hdots & 0 \\
             a_{21} & a_{22} & 0 &  \hdots & 0 \\
             \vdots & \vdots & \ddots & 0 & 0 \\
             a_{(n-1)1} & a_{(n-1)2} & \hdots & a_{(n-1)(n-1)} & 0 \\
             a_{n1} & a_{n2} & \hdots & a_{n(n-1)} & a_{nn}
             \end{bmatrix}\begin{pmatrix}
             x_1 \\
             x_2 \\
             \vdots\\
             x_{n-1} \\
             x_n
             \end{pmatrix}+\begin{pmatrix}
             x_1h_1(x_1,\ldots,x_n,t) \\
             h_2(x_1,\ldots,x_n,t) \\
             \vdots\\
             h_{n-1}(x_1,\ldots,x_n,t) \\
             h_n(x_1,\ldots,x_n,t)
             \end{pmatrix}\right\},\\
              (h_i \in {\Bbb C}(t)[x_1,\ldots,x_n], \ a_{ij} \in {\Bbb C}(t))
             \end{split}
             \end{align}
where $h_1$ is a polynomial which vanishes at $\overrightarrow{p}$ and $h_i$, $i=2,3,\ldots,n$ are polynomials of order at least 2 in $x_1,x_2,\ldots,x_n$, We call ordered set of the eigenvalues $(a_{11},a_{22},\cdots,a_{nn})$ {\it local index} at $\overrightarrow{p}$.

We are interested in the case with local index
\begin{equation}\label{integer}
(1,a_{22}/a_{11},\ldots,a_{nn}/a_{11}) \in {\Bbb Z}^{n}.
\end{equation}
These properties suggest the possibilities that $a_1$ is the residue of the formal Laurent series:
\begin{equation}
y_1(t)=\frac{a_{11}}{(t-t_0)}+b_1+b_2(t-t_0)+\cdots+b_n(t-t_0)^{n-1}+\cdots \quad (b_i \in {\Bbb C}),
\end{equation}
and the ratio $(1,a_{22}/a_{11},\ldots,a_{nn}/a_{11})$ is resonance data of the formal Laurent series of each $y_i(t) \ (i=2,\ldots,n)$, where $(y_1,\ldots,y_n)$ is original coordinate system satisfying $(x_1,\ldots,x_n)=(f_1(y_1,\ldots,y_n),\ldots,f_n(y_1,\ldots,y_n)), \ f_i(y_1,\ldots,y_n) \in {\Bbb C}(t)(y_1,\ldots,y_n)$.

If each component of $(1,a_{22}/a_{11},\ldots,a_{nn}/a_{11})$ has the same sign, we may resolve the accessible singularity by blowing-up finitely many times. However, when different signs appear, we may need to both blow up and blow down.

The $\alpha$-test,
\begin{equation}\label{poiuy}
t=t_0+\alpha T, \quad x_i=\alpha X_i, \quad \alpha \rightarrow 0,
\end{equation}
yields the following reduced system:
\begin{align}\label{ppppppp}
\begin{split}
\frac{d}{dT}\begin{pmatrix}
             X_1 \\
             X_2 \\
             \vdots\\
             X_{n-1} \\
             X_n
             \end{pmatrix}=\frac{1}{X_1}\begin{bmatrix}
             a_{11}(t_0) & 0 & 0 & \hdots & 0 \\
             a_{21}(t_0) & a_{22}(t_0) & 0 &  \hdots & 0 \\
             \vdots & \vdots & \ddots & 0 & 0 \\
             a_{(n-1)1}(t_0) & a_{(n-1)2}(t_0) & \hdots & a_{(n-1)(n-1)}(t_0) & 0 \\
             a_{n1}(t_0) & a_{n2}(t_0) & \hdots & a_{n(n-1)}(t_0) & a_{nn}(t_0)
             \end{bmatrix}\begin{pmatrix}
             X_1 \\
             X_2 \\
             \vdots\\
             X_{n-1} \\
             X_n
             \end{pmatrix},
             \end{split}
             \end{align}
where $a_{ij}(t_0) \in {\Bbb C}$. Fixing $t=t_0$, this system is the system of the first order ordinary differential equation with constant coefficient. Let us solve this system. At first, we solve the first equation:
\begin{equation}
X_1(T)=a_{11}(t_0)T+C_1 \quad (C_1 \in {\Bbb C}).
\end{equation}
Substituting this into the second equation in \eqref{ppppppp}, we can obtain the first order linear ordinary differential equation:
\begin{equation}
\frac{dX_2}{dT}=\frac{a_{22}(t_0) X_2}{a_{11}(t_0)T+C_1}+a_{21}(t_0).
\end{equation}
By variation of constant, in the case of $a_{11}(t_0) \not= a_{22}(t_0)$ we can solve explicitly:
\begin{equation}
X_2(T)=C_2(a_{11}(t_0)T+C_1)^{\frac{a_{22}(t_0)}{a_{11}(t_0)}}+\frac{a_{21}(t_0)(a_{11}(t_0)T+C_1)}{a_{11}(t_0)-a_{22}(t_0)} \quad (C_2 \in {\Bbb C}).
\end{equation}
This solution is a single-valued solution if and only if
$$
\frac{a_{22}(t_0)}{a_{11}(t_0)} \in {\Bbb Z}.
$$
In the case of $a_{11}(t_0)=a_{22}(t_0)$ we can solve explicitly:
\begin{equation}
X_2(T)=C_2(a_{11}(t_0)T+C_1)+\frac{a_{21}(t_0)(a_{11}(t_0)T+C_1){\rm Log}(a_{11}(t_0)T+C_1)}{a_{11}(t_0)} \quad (C_2 \in {\Bbb C}).
\end{equation}
This solution is a single-valued solution if and only if
$$
a_{21}(t_0)=0.
$$
Of course, $\frac{a_{22}(t_0)}{a_{11}(t_0)}=1 \in {\Bbb Z}$.
In the same way, we can obtain the solutions for each variables $(X_3,\ldots,X_n)$. The conditions $\frac{a_{jj}(t)}{a_{11}(t)} \in {\Bbb Z}, \ (j=2,3,\ldots,n)$ are necessary condition in order to have the Painlev\'e property.

\section{The case of the second Painlev\'e system}

In this section, we review the case of the second Painlev\'e system:
\begin{equation}\label{PII}
\frac{d^2u}{dt^2}=2u^3+tu+\alpha \quad (\alpha \in {\Bbb C}).
\end{equation}

Let us make its polynomial Hamiltonian from the viewpoint of accessible singularity and local index.

{\bf Step 0:} We make a change of variables.
\begin{equation}
x=u, \quad y=\frac{du}{dt}.
\end{equation}

{\bf Step 1:} We make a change of variables.
\begin{equation}
x_1=\frac{1}{x}, \quad y_1=\frac{y}{x^2}.
\end{equation}
In this coordinate system, we see that this system has two accessible singular points:
\begin{equation}
(x_1,y_1)=\left\{(0,1),(0,-1) \right\}.
\end{equation}

Around the point $(x_1,y_1)=(0,1)$, we can rewrite the system as follows.

{\bf Step 2:} We make a change of variables.
\begin{equation}
x_2=x_1, \quad y_2=y_1-1.
\end{equation}
In this coordinate system, we can rewrite the system satisfying the condition \eqref{b}:
\begin{align*}
\frac{d}{dt}\begin{pmatrix}
             x_2 \\
             y_2 
             \end{pmatrix}&=\frac{1}{x_2}\left\{\begin{pmatrix}
             -1 & 0   \\
             0 & -4 
             \end{pmatrix}\begin{pmatrix}
             x_2 \\
             y_2 
             \end{pmatrix}+\cdots\right\},
             \end{align*}
and we can obtain the local index $(-1,-4)$ at the point $\{(x_2,y_2)=(0,0)\}$. The ratio of the local index at the point $\{(x_2,y_2)=(0,0)\}$ is a positive integer.

We aim to obtain the local index $(-1,-2)$ by successive blowing-up procedures.

{\bf Step 3:} We blow up at the point $\{(x_2,y_2)=(0,0)\}$.
\begin{equation}
x_3=x_2, \quad y_3=\frac{y_2}{x_2}.
\end{equation}

{\bf Step 4:} We blow up at the point $\{(x_3,y_3)=(0,0)\}$.
\begin{equation}
x_4=x_3, \quad y_4=\frac{y_3}{x_3}.
\end{equation}
In this coordinate system, we see that this system has the following accessible singular point:
\begin{equation}
(x_4,y_4)=(0,t/2).
\end{equation}

{\bf Step 5:} We make a change of variables.
\begin{equation}
x_5=x_4, \quad y_5=y_4-t/2.
\end{equation}
In this coordinate system, we can rewrite the system as follows:
\begin{align*}
\frac{d}{dt}\begin{pmatrix}
             x_5 \\
             y_5 
             \end{pmatrix}&=\frac{1}{x_5}\left\{\begin{pmatrix}
             -1 & 0  \\
             \alpha-1/2 & -2
             \end{pmatrix}\begin{pmatrix}
             x_5 \\
             y_5 
             \end{pmatrix}+\cdots\right\},
             \end{align*}
and we can obtain the local index $(-1,-2)$. Here, the relation between $(x_5,y_5)$ and $(x,y)$ is given by
\begin{equation*}
  \left\{
  \begin{aligned}
   x_5 &=\frac{1}{x},\\
   y_5 &=y-x^2-\frac{t}{2}.
   \end{aligned}
  \right. 
\end{equation*}

Finally, we can choose canonical variables $(q,p)$.

{\bf Step 9:} We make a change of variables.
\begin{equation}
q=\frac{1}{x_5}, \quad p=y_5,
\end{equation}
and we can obtain the system
\begin{equation*}
  \left\{
  \begin{aligned}
   \frac{dq}{dt} &=q^2+p+\frac{t}{2},\\
   \frac{dp}{dt} &=-2qp+\alpha-\frac{1}{2}
   \end{aligned}
  \right. 
\end{equation*}
with the polynomial Hamiltonian $H_{II}$:
\begin{equation}
H_{II}=q^2 p+\frac{1}{2}p^2+\frac{t}{2}p-\left(\alpha-\frac{1}{2} \right)q.
\end{equation}
We remark that we can discuss the case of the accessible singular point $(x_1,y_1)=(0,-1)$ in the same way as in the case of $(x_1,y_1)=(0,1)$.

\section{Proof of theorem \ref{th1.1}}

By the same way of the second Painlev\'e system, we can prove Theorem \ref{th1.1}.

At first, we rewrite the equation \eqref{eq:1} to the system of the first order ordinary differential equations.

{\bf Step 0:} We make a change of variables.
\begin{equation}
x=u, \quad y=\frac{du}{dt}, \quad z=\frac{d^2u}{dt^2}, \quad w=\frac{d^3u}{dt^3}.
\end{equation}

{\bf Step 1:} We make a change of variables.
\begin{equation}
x_1=\frac{1}{x}, \quad y_1=\frac{y}{x^2}, \quad z_1=\frac{z}{x^3}, \quad w_1=\frac{w}{x^4}.
\end{equation}
In this coordinate system, we see that this system has four accessible singular points:
\begin{equation}
(x_1,y_1,z_1,w_1)=\left\{(0,-1,2,-6),\left(0,\frac{1}{2},\frac{1}{2},\frac{3}{4} \right),\left(0,\frac{1}{3},\frac{2}{9},\frac{2}{9} \right),\left(0,-\frac{1}{4},\frac{1}{8},-\frac{3}{32} \right) \right\}.
\end{equation}

Around the point $(x_1,y_1,z_1,w_1)=(0,-1,2,-6)$, we can rewrite the system as follows.
{\bf Step 2:} We make a change of variables.
\begin{equation}
x_2=x_1, \quad y_2=y_1+1, \quad z_2=z_1-2, \quad w_2=w_1+6.
\end{equation}
In this coordinate system, we can rewrite the system satisfying the condition \eqref{b}:
\begin{align*}
\frac{d}{dt}\begin{pmatrix}
             x_2 \\
             y_2 \\
             z_2 \\
             w_2
             \end{pmatrix}&=\frac{1}{x_2}\left\{\begin{pmatrix}
             1 & 0 & 0 & 0  \\
             0 & 4 & 1 & 0 \\
             0 & -6 & 3 & 1 \\
             0 & 4 & 10 & 4
             \end{pmatrix}\begin{pmatrix}
             x_2 \\
             y_2 \\
             z_2 \\
             w_2 
             \end{pmatrix}+\cdots\right\}.
             \end{align*}
To the above system, we make the linear transformation:
\begin{align*}
\begin{pmatrix}
             X_2 \\
             Y_2 \\
             Z_2 \\
             W_2
             \end{pmatrix}&=\begin{pmatrix}
             1 & 0 & 0 & 0  \\
             0 & 1 & 1 & 1 \\
             0 & -2 & -1 & 2 \\
             0 & 8 & 6 & 12
             \end{pmatrix}\begin{pmatrix}
             x_2 \\
             y_2 \\
             z_2 \\
             w_2 
             \end{pmatrix}
             \end{align*}
to arrive at
\begin{align*}
\frac{d}{dt}\begin{pmatrix}
             X_2 \\
             Y_2 \\
             Z_2 \\
             W_2
             \end{pmatrix}&=\frac{1}{X_2}\left\{\begin{pmatrix}
             1 & 0 & 0 & 0  \\
             0 & 2 & 0 & 0 \\
             0 & 0 & 3 & 0 \\
             0 & 0 & 0 & 6
             \end{pmatrix}\begin{pmatrix}
             X_2 \\
             Y_2 \\
             Z_2 \\
             W_2 
             \end{pmatrix}+\cdots\right\},
             \end{align*}
and we can obtain the local index $(1,2,3,6)$ at the point $\{(X_2,Y_2,Z_2,W_2)=(0,0,0,0)\}$. The continued ratio of the local index at the point $\{(x_2,y_2,z_2,w_2)=(0,0,0,0)\}$ are all positive integers
\begin{equation}
\left(\frac{2}{1},\frac{3}{1},\frac{6}{1} \right)=(2,3,6).
\end{equation}
This is the reason why we choose this accessible singular point.

We aim to obtain the local index $(1,0,0,2)$ by successive blowing-up procedures.

{\bf Step 3:} We blow up at the point $\{(x_2,y_2,z_2,w_2)=(0,0,0,0)\}$.
\begin{equation}
x_3=x_2, \quad y_3=\frac{y_2}{x_2}, \quad z_3=\frac{z_2}{x_2}, \quad w_3=\frac{w_2}{x_2}.
\end{equation}

{\bf Step 4:} We blow up at the point $\{(x_3,y_3,z_3,w_3)=(0,0,0,0)\}$.
\begin{equation}
x_4=x_3, \quad y_4=\frac{y_3}{x_3}, \quad z_4=\frac{z_3}{x_3}, \quad w_4=\frac{w_3}{x_3}.
\end{equation}
In this coordinate system, we see that this system has the following accessible singular locus:
\begin{equation}
(x_4,y_4,z_4,w_4)=(0,y_4,-2y_4,8y_4).
\end{equation}

{\bf Step 5:} We blow up along the curve $\{(x_4,y_4,z_4,w_4)=(0,y_4,-2y_4,8y_4)\}$.
\begin{equation}
x_5=x_4, \quad y_5=y_4, \quad z_5=\frac{z_4+2y_4}{x_4}, \quad w_5=\frac{w_4-8y_4}{x_4}.
\end{equation}
In this coordinate system, we see that this system has the following accessible singular locus:
\begin{equation}
(x_5,y_5,z_5,w_5)=(0,y_5,z_5,-2z_5).
\end{equation}

{\bf Step 6:} We blow up along the surface $\{(x_5,y_5,z_5,w_5)=(0,y_5,z_5,-2z_5)\}$.
\begin{equation}
x_6=x_5, \quad y_6=y_5, \quad z_6=z_5, \quad w_6=\frac{w_5+2z_5}{x_5}.
\end{equation}
In this coordinate system, we see that this system has the following accessible singular locus:
\begin{equation}
(x_6,y_6,z_6,w_6)=\left(0,y_6,z_6,-\frac{3}{2}y_6^2-\frac{t}{2} \right).
\end{equation}

{\bf Step 7:} We make a change of variables.
\begin{equation}
x_7=x_6, \quad y_7=y_6, \quad z_7=z_6, \quad w_7=w_6+\frac{3}{2}y_6^2+\frac{t}{2}.
\end{equation}
In this coordinate system, we can rewrite the system as follows:
\begin{align*}
\frac{d}{dt}\begin{pmatrix}
             x_7 \\
             y_7 \\
             z_7 \\
             w_7
             \end{pmatrix}&=\frac{1}{x_7}\left\{\begin{pmatrix}
             1 & 0 & 0 & 0  \\
             0 & 0 & 0 & 0 \\
             -\frac{t}{2} & 0 & 0 & 0 \\
             \alpha+\frac{1}{2} & 0 & 0 & 2
             \end{pmatrix}\begin{pmatrix}
             x_7 \\
             y_7 \\
             z_7 \\
             w_7 
             \end{pmatrix}+\cdots\right\},
             \end{align*}
and we can obtain the local index $(1,0,0,2)$. Here, the relation between $(x_7,y_7,z_7,w_7)$ and $(x,y,z,w)$ is given by
\begin{equation*}
  \left\{
  \begin{aligned}
   x_7 &=\frac{1}{x},\\
   y_7 &=x^2+y,\\
   z_7 &=z+2xy,\\
   w_7 &=w+\frac{t}{2}-\frac{1}{2}x^4-x^2 y+\frac{3}{2}y^2+2xz.
   \end{aligned}
  \right. 
\end{equation*}

{\bf Step 8:} We make a change of variables.
\begin{equation}
x_8=\frac{1}{x_7}, \quad y_8=y_7, \quad z_8=z_7, \quad w_8=w_7.
\end{equation}
In this coordinate system, we can rewrite the system as follows:
\begin{equation*}
  \left\{
  \begin{aligned}
   \frac{dx_8}{dt} &=-x_8^2+y_8,\\
    \frac{dy_8}{dt} &=z_8,\\
    \frac{dz_8}{dt} &=\frac{1}{2}y_8^2+w_8-\frac{t}{2},\\
    \frac{dw_8}{dt} &=2x_8 w_8+\alpha+\frac{1}{2}.
   \end{aligned}
  \right. 
\end{equation*}
Finally, we can choose canonical variables $(q_1,p_1,q_2,p_2)$.

{\bf Step 9:} We make a change of variables.
\begin{equation}
q_1=-x_8, \quad p_1=-w_8, \quad q_2=-z_8, \quad p_2=-y_8,
\end{equation}
and we can obtain the system \eqref{eq:3} with the polynomial Hamiltonian \eqref{eq:4}.

Thus, we have completed the proof of Theorem \ref{th1.1} \qed.

We note on the remaining accessible singular points.

Around the point $(x_1,y_1,z_1,w_1)=(0,1/2,1/2,3/4)$, we can rewrite the system as follows.

{\bf Step 2:} We make a change of variables.
\begin{equation}
x_2=x_1, \quad y_2=y_1-1/2, \quad z_2=z_1-1/2, \quad w_2=w_1-3/4.
\end{equation}
In this coordinate system, we can rewrite the system satisfying the condition \eqref{b}:
\begin{align*}
\frac{d}{dt}\begin{pmatrix}
             x_2 \\
             y_2 \\
             z_2 \\
             w_2
             \end{pmatrix}&=\frac{1}{x_2}\left\{\begin{pmatrix}
             -1/2 & 0 & 0 & 0  \\
             0 & -2 & 1 & 0 \\
             0 & -3/2 & -3/2 & 1 \\
             0 & -1/2 & 5/2 & -2
             \end{pmatrix}\begin{pmatrix}
             x_2 \\
             y_2 \\
             z_2 \\
             w_2 
             \end{pmatrix}+\cdots\right\}.
             \end{align*}
To the above system, we make the linear transformation:
\begin{align*}
\begin{pmatrix}
             X_2 \\
             Y_2 \\
             Z_2 \\
             W_2
             \end{pmatrix}&=\begin{pmatrix}
             1 & 0 & 0 & 0  \\
             0 & 1/2 & 2/3 & 1/3 \\
             0 & 1/2 & 1/3 & -1/3 \\
             0 & 1 & 1 & 1
             \end{pmatrix}\begin{pmatrix}
             x_2 \\
             y_2 \\
             z_2 \\
             w_2 
             \end{pmatrix}
             \end{align*}
to arrive at
\begin{align*}
\frac{d}{dt}\begin{pmatrix}
             X_2 \\
             Y_2 \\
             Z_2 \\
             W_2
             \end{pmatrix}&=\frac{1}{X_2}\left\{\begin{pmatrix}
             -1/2 & 0 & 0 & 0  \\
             0 & -1 & 0 & 0 \\
             0 & 0 & -3/2 & 0 \\
             0 & 0 & 0 & -3
             \end{pmatrix}\begin{pmatrix}
             X_2 \\
             Y_2 \\
             Z_2 \\
             W_2 
             \end{pmatrix}+\cdots\right\},
             \end{align*}
and we can obtain the local index $(-1/2,-1,-3/2,-3)$ at the point $\{(X_2,Y_2,Z_2,W_2)=(0,0,0,0)\}$. The continued ratio of the local index at the point $\{(x_2,y_2,z_2,w_2)=(0,0,0,0)\}$ are all positive integers
\begin{equation}
\left(\frac{-1}{-1/2},\frac{-3/2}{-1/2},\frac{-3}{-1/2} \right)=(2,3,6).
\end{equation}
We remark that we can discuss this case in the same way as in the case of $(x_1,y_1,z_1,w_1)=(0,-1,2,-6)$.

Around the point $(x_1,y_1,z_1,w_1)=(0,1/3,2/9,2/9)$, we can rewrite the system as follows.

{\bf Step 2:} We make a change of variables.
\begin{equation}
x_2=x_1, \quad y_2=y_1-1/3, \quad z_2=z_1-2/9, \quad w_2=w_1-2/9.
\end{equation}
In this coordinate system, we can rewrite the system satisfying the condition \eqref{b}:
\begin{align*}
\frac{d}{dt}\begin{pmatrix}
             x_2 \\
             y_2 \\
             z_2 \\
             w_2
             \end{pmatrix}&=\frac{1}{x_2}\left\{\begin{pmatrix}
             -1/3 & 0 & 0 & 0  \\
             0 & -4/3 & 1 & 0 \\
             0 & -2/3 & -1 & 1 \\
             0 & 4/3 & 10/3 & -4/3
             \end{pmatrix}\begin{pmatrix}
             x_2 \\
             y_2 \\
             z_2 \\
             w_2 
             \end{pmatrix}+\cdots\right\}.
             \end{align*}
To the above system, we make the linear transformation:
\begin{align*}
\begin{pmatrix}
             X_2 \\
             Y_2 \\
             Z_2 \\
             W_2
             \end{pmatrix}&=\begin{pmatrix}
             1 & 0 & 0 & 0  \\
             0 & 1/4 & 3/4 & 1/2 \\
             0 & 1/2 & -1/2 & -1/2 \\
             0 & 1 & 1 & 1
             \end{pmatrix}\begin{pmatrix}
             x_2 \\
             y_2 \\
             z_2 \\
             w_2 
             \end{pmatrix}
             \end{align*}
to arrive at
\begin{align*}
\frac{d}{dt}\begin{pmatrix}
             X_2 \\
             Y_2 \\
             Z_2 \\
             W_2
             \end{pmatrix}&=\frac{1}{X_2}\left\{\begin{pmatrix}
             -1/3 & 0 & 0 & 0  \\
             0 & 2/3 & 0 & 0 \\
             0 & 0 & -2 & 0 \\
             0 & 0 & 0 & -7/3
             \end{pmatrix}\begin{pmatrix}
             X_2 \\
             Y_2 \\
             Z_2 \\
             W_2 
             \end{pmatrix}+\cdots\right\},
             \end{align*}
and we can obtain the local index $(-1/3,2/3,-2,-7/3)$ at the point $\{(X_2,Y_2,Z_2,W_2)=(0,0,0,0)\}$. The continued ratio of the local index at the point $\{(x_2,y_2,z_2,w_2)=(0,0,0,0)\}$ are
\begin{equation}
\left(\frac{2/3}{-1/3},\frac{-2}{-1/3},\frac{-7/3}{-1/3} \right)=(-2,6,7).
\end{equation}
In this case, the local index involves a negative integer. So, we need to blow down.

Around the point $(x_1,y_1,z_1,w_1)=(0,-1/4,1/8,-3/32)$, we can rewrite the system as follows.

{\bf Step 2:} We make a change of variables.
\begin{equation}
x_2=x_1, \quad y_2=y_1+1/4, \quad z_2=z_1-1/8, \quad w_2=w_1+3/32.
\end{equation}
In this coordinate system, we can rewrite the system satisfying the condition \eqref{b}:
\begin{align*}
\frac{d}{dt}\begin{pmatrix}
             x_2 \\
             y_2 \\
             z_2 \\
             w_2
             \end{pmatrix}&=\frac{1}{x_2}\left\{\begin{pmatrix}
             1/4 & 0 & 0 & 0  \\
             0 & 1 & 1 & 0 \\
             0 & -3/8 & 3/4 & 1 \\
             0 & -11/4 & 25/4 & 1
             \end{pmatrix}\begin{pmatrix}
             x_2 \\
             y_2 \\
             z_2 \\
             w_2 
             \end{pmatrix}+\cdots\right\}.
             \end{align*}
To the above system, we make the linear transformation:
\begin{align*}
\begin{pmatrix}
             X_2 \\
             Y_2 \\
             Z_2 \\
             W_2
             \end{pmatrix}&=\begin{pmatrix}
             1 & 0 & 0 & 0  \\
             0 & 8/39 & 4/29 & 4/3 \\
             0 & 16/39 & -11/29 & 2/3 \\
             0 & 1 & 1 & 1
             \end{pmatrix}\begin{pmatrix}
             x_2 \\
             y_2 \\
             z_2 \\
             w_2 
             \end{pmatrix}
             \end{align*}
to arrive at
\begin{align*}
\frac{d}{dt}\begin{pmatrix}
             X_2 \\
             Y_2 \\
             Z_2 \\
             W_2
             \end{pmatrix}&=\frac{1}{X_2}\left\{\begin{pmatrix}
             1/4 & 0 & 0 & 0  \\
             0 & 3 & 0 & 0 \\
             0 & 0 & -7/4 & 0 \\
             0 & 0 & 0 & 3/2
             \end{pmatrix}\begin{pmatrix}
             X_2 \\
             Y_2 \\
             Z_2 \\
             W_2 
             \end{pmatrix}+\cdots\right\},
             \end{align*}
and we can obtain the local index $(1/4,3,-7/4,3/2)$ at the point $\{(X_2,Y_2,Z_2,W_2)=(0,0,0,0)\}$. The continued ratio of the local index at the point $\{(x_2,y_2,z_2,w_2)=(0,0,0,0)\}$ are
\begin{equation}
\left(\frac{3}{1/4},\frac{-7/4}{1/4},\frac{3/2}{1/4} \right)=(12,-7,6).
\end{equation}
In this case, the local index involves a negative integer. So, we need to blow down.

\section{Symmetry and holomorphy conditions}
In this section, we study the symmetry and holomorphy conditions of the system \eqref{eq:3}. These symmetries, holomorphy conditions and invariant divisors are new.
\begin{theorem}\label{pro:3}
Let us consider a polynomial Hamiltonian system with Hamiltonian $H \in {\Bbb C}(t)[q_1,p_1,q_2,p_2]$. We assume that

$(A1)$ $deg(H)=5$ with respect to $q_1,p_1,q_2,p_2$.

$(A2)$ This system becomes again a polynomial Hamiltonian system in each coordinate system $r_i \ (i=0,1)${\rm : \rm}
\begin{align}
\begin{split}
r_0:&x_0=\frac{1}{q_1}, \ y_0=-\left(f_0p_1+\left(\alpha+\frac{1}{2} \right) \right)q_1, \ z_0=q_2, \ w_0=p_2, \\
r_1:&x_1=\frac{1}{q_1}, \ y_1=-\left(f_1q_1-2(\alpha-1) \right)q_1, \ z_1=q_2+3q_1^3+3q_1p_2, \ w_1=p_2+\frac{3}{2} q_1^2,
\end{split}
\end{align}
where $f_0:=p_1$ and $f_1:=p_1+3q_1q_2-\frac{3}{2}\left(p_2^2-t \right)$. Then such a system coincides with the system \eqref{eq:3}  with the polynomial Hamiltonian \eqref{eq:4}.
\end{theorem}
We note that the condition $(A2)$ should be read that
\begin{align*}
&r_0(H), \quad r_1 \left(H-\frac{3}{2} q_1 \right)
\end{align*}
are polynomials with respect to $x_i,y_i,z_i,w_i$.

\begin{theorem}\label{th:2}
The system \eqref{eq:3} is invariant under the following transformations $s_0,s_1:$ with {\it the notation} $(*):=(q_1,p_1,q_2,p_2,t;\alpha)$\rm{; \rm}
\begin{align}
\begin{split}
s_0:(*) \rightarrow &\left(q_1+\frac{\alpha+\frac{1}{2}}{f_0},p_1,q_2,p_2,t;-1-\alpha \right),\\
s_1:(*) \rightarrow &(q_1+\frac{2-2\alpha}{f_1}, \quad p_1+\frac{6q_2(\alpha-1)}{f_1}+\frac{18p_2(\alpha-1)^2}{f_1^2}+\frac{36q_1(\alpha-1)^3}{f_1^3}\\
&-\frac{18(\alpha-1)^4}{f_1^4}, \quad q_2+\frac{6p_2(\alpha-1)}{f_1}+\frac{18q_1(\alpha-1)^2}{f_1^2}-\frac{12(\alpha-1)^3}{f_1^3}\\
&p_2+\frac{6q_1(\alpha-1)}{f_1}-\frac{6(\alpha-1)^2}{f_1^2},t;2-\alpha).
\end{split}
\end{align}
\end{theorem}
The B{\"a}cklund transformations of the system \eqref{eq:3} satisfy
\begin{equation}
s_i(g)=g+\frac{\alpha_i}{f_i}\{f_i,g\}+\frac{1}{2!} \left(\frac{\alpha_i}{f_i} \right)^2 \{f_i,\{f_i,g\} \}+\cdots \quad (g \in {\Bbb C}(t)[q_1,p_1,q_2,p_2]),
\end{equation}
where poisson bracket $\{,\}$ satisfies the relations:
$$
\{p_1,q_1\}=\{p_2,q_2\}=1, \quad the \ others \ are \ 0.
$$
Since these B{\"a}cklund transformations have Lie theoretic origin, similarity reduction of a Drinfeld-Sokolov hierarchy admits such a B{\"a}cklund symmetry.

In the next section, we introduce a symmetric form in dimension five for this Hamiltonian system by taking the two invariant divisors $f_0,f_1$ as the dependent variables.

\section{Symmetric form}

In this section, we find a 1-parameter family of  ordinary differential systems in dimension five with affine Weyl group symmetry of type $A_2^{(2)}$ explicitly given by
\begin{equation}\label{eq:11}
  \left\{
  \begin{aligned}
\frac{dx}{dt}=&-2xz-\frac{3}{2}\alpha_0,\\
\frac{dy}{dt}=&yz+\frac{3}{2}\alpha_1,\\
\frac{dz}{dt}=&z^2+q,\\
\frac{dw}{dt}=&-zw+\frac{2}{3}x+\frac{1}{3}y,\\
\frac{dq}{dt}=&w.
   \end{aligned}
  \right. 
\end{equation}
Here $x,y,z,w$ and $q$ denote unknown complex variables, and $\alpha_0,\alpha_1$ are complex parameters satisfying the relation:
\begin{equation}\label{eq:12}
2\alpha_0+\alpha_1=1.
\end{equation}
This system can be considered as giving a symmetric form for the Hamiltonian system \eqref{eq:3} by taking the two invariant divisors $x,y$ as the dependent variables (see \eqref{inv}).

\begin{theorem}\label{th:1}
The system \eqref{eq:11} admits the affine Weyl group symmetry of type $A_2^{(2)}$ as the group of its B{\"a}cklund transformations, whose generators $s_0,s_1$ defined as follows$:$ with {\it the notation} $(*):=(x,y,z,w,q;\alpha_0,\alpha_1)$,
\begin{align*}
s_0:(*) \rightarrow &\left(x,y+\frac{9\alpha_0 w}{x},z+\frac{3\alpha_0}{x},w,q;-\alpha_0,\alpha_1+4\alpha_0 \right),\\
s_1:(*) \rightarrow &(x-\frac{9\alpha_1 w}{y}+\frac{162\alpha_1^2 q}{4y^2}-\frac{243\alpha_1^3 z}{2y^3}-\frac{729\alpha_1^4}{8y^4},y,z+\frac{3\alpha_1}{y},\\
&w-\frac{9\alpha_1 q}{y}+\frac{81\alpha_1^2 z}{2y^2}+\frac{81\alpha_1^3}{2y^3},q-\frac{9\alpha_1 z}{y}-\frac{27\alpha_1^2}{2y^2};\\
&\alpha_0+\alpha_1,-\alpha_1).
\end{align*}
\end{theorem}

\begin{proposition}
The system \eqref{eq:11} has the following invariant divisors\rm{:\rm}
\begin{center}\label{inv}
\begin{tabular}{|c|c|} \hline
parameter's relation  & invariant divisor \\ \hline
$\alpha_0=0$ & $x$ \\ \hline
$\alpha_1=0$ & $y$ \\ \hline
\end{tabular}
\end{center}
\end{proposition}
We note that when $\alpha_0=0$, we see that the system \eqref{eq:11} admits a particular solution $x=0$.

\begin{theorem}
Let us consider the following ordinary differential system in the polynomial class\rm{:\rm}
\begin{align*}
&\frac{dx}{dt}=f_1(x,y,z,w,q), \cdots, \frac{dq}{dt}=f_5(x,y,z,w,q) \quad (f_i \in {\Bbb C}(t)[x,y,z,w,q]).
\end{align*}
We assume that

$(A1)$ $deg(f_i)=3$ with respect to $x,y,z,w,q$.

$(A2)$ The right-hand side of this system becomes again a polynomial in each coordinate system $(x_i,y_i,z_i,w_i,q_i) \ (i=0,1):$

\begin{align}
\begin{split}
0) \ &x_0=-\left(xz+3\alpha_0 \right)z, \quad y_0=\frac{y}{z}, \quad z_0=\frac{1}{z}, \quad w_0=\left(w-\frac{y}{3z} \right)z, \quad q_0=q,\\
1) \ &x_1=x+3zw+\frac{9}{2}z^2 q+\frac{27}{8}z^4, \quad y_1=-(yz+3\alpha_1)z, \quad z_1=\frac{1}{z},\\
&w_1=w+3zq+3z^3, \quad q_1=q+\frac{3}{2}z^2.
\end{split}
\end{align}
Then such a system coincides with the system \eqref{eq:11}.
\end{theorem}
We note that these transition functions satisfy the condition{\rm:\rm}
\begin{align*}
&dx_i \wedge dy_i \wedge dz_i \wedge dw_i \wedge dq_i=dx \wedge dy \wedge dz \wedge dw \wedge dq \quad (i=0,1).
\end{align*}

For this system let us try to seek its first integrals of polynomial type with respect to $x,y,z,w,q$.
\begin{proposition}
This system \eqref{eq:11} has its first integral\rm{:\rm}
\begin{equation*}
y-x-3zw+\frac{3}{2}q^2-\frac{3}{2}t=0.
\end{equation*}
\end{proposition}

By using this, elimination of $y$ from the system \eqref{eq:11} coincides with the system \eqref{eq:3}.
\begin{proposition}
By the transformations
\begin{equation*}
q_1 =z, \quad p_1 =x, \quad q_2 =w, \quad p_2 =q, \quad \alpha=\frac{1}{2}(6\alpha_0-1),
\end{equation*}
elimination of $y$ from the system \eqref{eq:11} coincides with the system \eqref{eq:3} with the polynomial Hamiltonian \eqref{eq:4}.
\end{proposition}

\section{Autonomous version of the system \eqref{eq:3} and partial differential system in two variables}

In this section, we find an autonomous version of the system  \eqref{eq:3} given by
\begin{equation}\label{eq:A1}
  \left\{
  \begin{aligned}
   dq_1 =&\frac{\partial K_1}{\partial p_1}dt+\frac{\partial K_2}{\partial p_1}ds,\\
   dp_1 =&-\frac{\partial K_1}{\partial q_1}dt-\frac{\partial K_2}{\partial q_1}ds,\\
   dq_2 =&\frac{\partial K_1}{\partial p_2}dt+\frac{\partial K_2}{\partial p_2}ds,\\
   dp_2 =&-\frac{\partial K_1}{\partial q_2}dt-\frac{\partial K_2}{\partial q_2}ds
   \end{aligned}
  \right. 
\end{equation}
with the polynomial Hamiltonians
\begin{align}\label{eq:A2}
\begin{split}
K_1=&q_1^2p_1+3\alpha_0 q_1-\frac{q_2^2}{2}-\frac{p_2^3}{6}+p_1 p_2,\\
K_2=&\frac{p_1^3}{9}-\frac{3\alpha_0}{2}q_2 p_2^2+3\alpha_0^2 p_2+\alpha_0 p_1q_2+3\alpha_0 q_1q_2^2-\frac{1}{3}p_1^2 p_2^2+\frac{2}{3}q_1p_1^2q_2+\frac{1}{4}p_1p_2^4\\
&+q_1^2p_1q_2^2-q_1p_1q_2p_2^2.
\end{split}
\end{align}

\begin{proposition}
The system \eqref{eq:A1} satisfies the compatibility conditions$:$
\begin{equation}
\frac{\partial }{\partial s} \frac{\partial q_1}{\partial t}=\frac{\partial }{\partial t} \frac{\partial q_1}{\partial s}, \quad \frac{\partial }{\partial s} \frac{\partial p_1}{\partial t}=\frac{\partial }{\partial t} \frac{\partial p_1}{\partial s}, \quad \frac{\partial }{\partial s} \frac{\partial q_2}{\partial t}=\frac{\partial }{\partial t} \frac{\partial q_2}{\partial s}, \quad \frac{\partial }{\partial s} \frac{\partial p_2}{\partial t}=\frac{\partial }{\partial t} \frac{\partial p_2}{\partial s}.
\end{equation}
\end{proposition}

\begin{proposition}
The system \eqref{eq:A1} has $K_1$ and $K_2$ as its first integrals.
\end{proposition}

\begin{proposition}
Two Hamiltonians $K_1$ and $K_2$ satisfy
\begin{equation}
\{K_1,K_2\}=0,
\end{equation}
where
\begin{equation}
\{K_1,K_2\}=\frac{\partial K_1}{\partial p_1}\frac{\partial K_2}{\partial q_1}-\frac{\partial K_1}{\partial q_1}\frac{\partial K_2}{\partial p_1}+\frac{\partial K_1}{\partial p_2}\frac{\partial K_2}{\partial q_2}-\frac{\partial K_1}{\partial q_2}\frac{\partial K_2}{\partial p_2}.
\end{equation}
\end{proposition}
Here, $\{,\}$ denotes the poisson bracket such that $\{p_i,q_j\}={\delta}_{ij}$ (${\delta}_{ij}$:kronecker's delta).

\begin{theorem}
The system \eqref{eq:A1} admits the affine Weyl group symmetry of type $A_2^{(2)}$ as the group of its B{\"a}cklund transformations, whose generators $s_0,s_1$ defined as follows$:$ with {\it the notation} $(*):=(q_1,p_1,q_2,p_2,t,s;\alpha_0,\alpha_1)$\rm{; \rm}
\begin{align*}
s_0:(*) \rightarrow &\left(q_1+\frac{3\alpha_0}{p_1},p_1,q_2,p_2,t,s;-\alpha_0,\alpha_1+4\alpha_0 \right),\\
s_1:(*) \rightarrow &(q_1-\frac{6\alpha_1}{3p_2^2-2p_1-6q_1q_2},\\
&p_1+\frac{18\alpha_1q_2}{3p_2^2-2p_1-6q_1q_2}+\frac{162\alpha_1^2p_2}{(3p_2^2-2p_1-6q_1q_2)^2}+\frac{972\alpha_1^3q_1}{(3p_2^2-2p_1-6q_1q_2)^3}\\
&-\frac{1458\alpha_1^4}{(3p_2^2-2p_1-6q_1q_2)^4},\\
&q_2+\frac{18\alpha_1p_2}{3p_2^2-2p_1-6q_1q_2}+\frac{162\alpha_1^2 q_1}{(3p_2^2-2p_1-6q_1q_2)^2}-\frac{324\alpha_1^3}{(3p_2^2-2p_1-6q_1q_2)^3},\\
&p_2+\frac{18\alpha_1 q_1}{3p_2^2-2p_1-6q_1q_2}-\frac{54\alpha_1^2}{(3p_2^2-2p_1-6q_1q_2)^2},t,s;\alpha_0+\alpha_1,-\alpha_1).
\end{align*}
\end{theorem}
Here, the parameters $\alpha_i$ satisfy the relation $2\alpha_0+\alpha_1=0$.

\begin{theorem}
Let us consider a polynomial Hamiltonian system with Hamiltonian $K \in {\Bbb C}[q_1,p_1,q_2,p_2]$. We assume that

$(C1)$ $deg(K)=5$ with respect to $q_1,p_1,q_2,p_2$.

$(C2)$ This system becomes again a polynomial Hamiltonian system in each coordinate $R_i \ (i=0,1)${\rm : \rm}
\begin{align*}
\begin{split}
R_0:(x_0,y_0,z_0,w_0)=&\left(\frac{1}{q_1},-(q_1f_0+3\alpha_0)q_1,q_2,p_2 \right),\\
R_1:(x_1,y_1,z_1,w_1)=&\left(\frac{1}{q_1},-(q_1f_1+3\alpha_1)q_1,q_2+3q_1p_2+3q_1^3,p_2+\frac{3}{2}q_1^2 \right),
\end{split}
\end{align*}
where $f_0:=p_1$ and $f_1:=p_1-\frac{3}{2}p_2^2+3q_1q_2$. Then such a system coincides with the Hamiltonian system \eqref{eq:A1} with the polynomial Hamiltonians $K_1,K_2$.
\end{theorem}
We note that the conditions $(C2)$ should be read that
\begin{align*}
\begin{split}
&R_0(K), \quad R_1(K)
\end{split}
\end{align*}
are polynomials with respect to $x_i,y_i,z_i,w_i$.

Next, let us consider the relation between the polynomial Hamiltonian system \eqref{eq:A1} and autonomous version of the system  \eqref{eq:3}. In this paper, we can make the birational transformations between the polynomial Hamiltonian system \eqref{eq:A1} and autonomous version of the system  \eqref{eq:3}.
\begin{theorem}
The birational transformations
\begin{equation}\label{eq:A3}
  \left\{
  \begin{aligned}
   x =&-q_1,\\
   y =&-p_2-q_1^2,\\
   z =&-q_2-2q_1p_2-2q_1^3,\\
   w =&-p_1-2q_1q_2-\frac{3}{2}p_2^2-8q_1^2p_2-6q_1^4
   \end{aligned}
  \right. 
\end{equation}
take the Hamiltonian system \eqref{eq:A1} to the system
\begin{equation}\label{eq:A4}
  \left\{
  \begin{aligned}
   dx =&y dt+(-\frac{w^2}{3}-\frac{x^4 w}{3}+3\alpha_0 xy+2x^2yw+\frac{2x^6 y}{3}-\frac{5y^2 w}{3}-\frac{10x^4 y^2}{3}+\frac{14x^2y^3}{3}-2y^4\\
&+\alpha_0 z-\frac{x^5z}{3}+\frac{2x^3yz}{3}-\frac{xy^2z}{3}+\frac{x^2z^2}{3})ds,\\
   dy =&z dt+f_1(x,y,z,w)ds,\\
   dz =&w dt+f_2(x,y,z,w)ds,\\
   dw =&\left(5x^2 z-x^5+5xy^2-5yz+3\alpha_0 \right)dt+f_3(x,y,z,w)ds,
   \end{aligned}
  \right. 
\end{equation}
where $f_i(x,y,z,w) \in {\Bbb C}[x,y,z,w] \ (i=1,2,3)$. 
\end{theorem}
Setting $u:=x$, we see that
\begin{equation}
\frac{\partial u}{\partial t}=y, \quad  \frac{\partial^2 u}{\partial t^2}=z, \quad \frac{\partial^3 u}{\partial t^3}=w,
\end{equation}
and
\begin{equation}\label{eq:A5}
  \left\{
  \begin{aligned}
   \frac{\partial^4 u}{\partial t^4} =&5u^2 \frac{\partial^2 u}{\partial t^2}-u^5+5u\left(\frac{\partial u}{\partial t} \right)^2-5\frac{\partial u}{\partial t}\frac{\partial^2 u}{\partial t^2}+3\alpha_0,\\
   \frac{\partial u}{\partial s} =&-\frac{1}{3}\left(\frac{\partial^3 u}{\partial t^3} \right)^2-\frac{1}{3}u^4\frac{\partial^3 u}{\partial t^3} +3\alpha_0 u\frac{\partial u}{\partial t}+2u^2\frac{\partial u}{\partial t}\frac{\partial^3 u}{\partial t^3}+\frac{2}{3}u^6 \frac{\partial u}{\partial t}-\frac{5}{3}\left(\frac{\partial u}{\partial t} \right)^2 \frac{\partial^3 u}{\partial t^3}\\
&-\frac{10}{3}u^4 \left(\frac{\partial u}{\partial t} \right)^2+\frac{14}{3}u^2\left(\frac{\partial u}{\partial t} \right)^3-2\left(\frac{\partial u}{\partial t} \right)^4+\alpha_0 \frac{\partial^2 u}{\partial t^2}-\frac{1}{3}u^5 \frac{\partial^2 u}{\partial t^2}+\frac{2}{3}u^3\frac{\partial u}{\partial t}\frac{\partial^2 u}{\partial t^2}\\
&-\frac{1}{3}u\left(\frac{\partial u}{\partial t} \right)^2\frac{\partial^2 u}{\partial t^2}+\frac{1}{3}u^2\left(\frac{\partial^2 u}{\partial t^2} \right)^2.
   \end{aligned}
  \right. 
\end{equation}
The first equation in \eqref{eq:A5} coincides with an autonomous version of the system  \eqref{eq:3}.

We see that the second equation in \eqref{eq:A5} can be considered as homogeneous polynomial of degree 8 when setting $[\alpha_0]=5$ and $[s]=7$.

Making a partial derivation in the variable $t$ for the first equation of \eqref{eq:A5}, we obtain
\begin{equation}\label{eq:A6}
  \left\{
  \begin{aligned}
   \frac{\partial^5 u}{\partial t^5} =&-5\left(\frac{\partial u}{\partial t}-u^2 \right)\frac{\partial^3 u}{\partial t^3}-5\left(\frac{\partial^2 u}{\partial t^2} \right)^2+20u\frac{\partial u}{\partial t}\frac{\partial^2 u}{\partial t^2}+5\left(\frac{\partial u}{\partial t} \right)^3-5u^4\frac{\partial u}{\partial t},\\
   \frac{\partial u}{\partial s} =&-\frac{1}{3}\left(\frac{\partial^3 u}{\partial t^3} \right)^2-\frac{1}{3}u^4\frac{\partial^3 u}{\partial t^3} +3\alpha_0 u\frac{\partial u}{\partial t}+2u^2\frac{\partial u}{\partial t}\frac{\partial^3 u}{\partial t^3}+\frac{2}{3}u^6 \frac{\partial u}{\partial t}-\frac{5}{3}\left(\frac{\partial u}{\partial t} \right)^2 \frac{\partial^3 u}{\partial t^3}\\
&-\frac{10}{3}u^4 \left(\frac{\partial u}{\partial t} \right)^2+\frac{14}{3}u^2\left(\frac{\partial u}{\partial t} \right)^3-2\left(\frac{\partial u}{\partial t} \right)^4+\alpha_0 \frac{\partial^2 u}{\partial t^2}-\frac{1}{3}u^5 \frac{\partial^2 u}{\partial t^2}+\frac{2}{3}u^3\frac{\partial u}{\partial t}\frac{\partial^2 u}{\partial t^2}\\
&-\frac{1}{3}u\left(\frac{\partial u}{\partial t} \right)^2\frac{\partial^2 u}{\partial t^2}+\frac{1}{3}u^2\left(\frac{\partial^2 u}{\partial t^2} \right)^2.
   \end{aligned}
  \right. 
\end{equation}
The first equation in \eqref{eq:A6} coincides with the Modified Sawada-Kotera equation \eqref{SKeq} satisfying the relation
$$
u:=g, \quad \frac{\partial g}{\partial s}=0, \quad w:=t.
$$

\begin{question}
It is still an open question whether the second equation in \eqref{eq:A6} coincides with which of the soliton equations.
\end{question}

\end{document}